\documentclass[11pt,a4paper]{amsart}
\usepackage[all]{xy}

\newcommand{\ie}{{\em i.e.}\ }
\newcommand{\cf}{{\em cf.}\ }
\newcommand{\eg}{{\em e.g.}\ }
\newcommand{\ko}{\: , \;}

\renewcommand{\tilde}[1]{\widetilde{#1}}

\newtheorem*{theorem}{Theorem}
\newtheorem*{lemma}{Lemma}
\newtheorem*{proposition}{Proposition}

\newcommand{\opname}[1]{\operatorname{\mathsf{#1}}}

\renewcommand{\mod}{\opname{mod}\nolimits}

\newcommand{\ind}{\opname{ind}\nolimits}
\newcommand{\indec}{\opname{indec}\nolimits}

\newcommand{\Z}{\mathbb{Z}}
\newcommand{\Q}{\mathbb{Q}}

\newcommand{\iso}{\stackrel{_\sim}{\rightarrow}}
\newcommand{\id}{\mathbf{1}}

\newcommand{\Hom}{\opname{Hom}}

\newcommand{\Ext}{\opname{Ext}}

\newcommand{\Aut}{\opname{Aut}}
\newcommand{\Lie}{\opname{Lie}}

\newcommand{\ca}{{\mathcal A}}

\newcommand{\cc}{{\mathcal C}}

\newcommand{\cF}{{\mathcal F}}

\newcommand{\ct}{{\mathcal T}}

\newcommand{\eps}{\varepsilon}

\newcommand{\sgn}{\mbox{sgn}}

\begin{document}


\title{On the combinatorics of rigid objects in $2$-Calabi-Yau categories}
\author{Raika Dehy and Bernhard Keller}

\address{Universit\'e Cergy-Pontoise/Saint-Martin,
D\'epartment de Math\'ematiques, UMR 8088 du CNRS, 2 avenue Adolphe
Chauvin, 95302 Cergy-Pontoise Cedex, France}
\email{Raika.Dehy@math.u-cergy.fr}

\address{UFR de Math\'ematiques, Universit\'e Denis Diderot --
Paris 7, Institut de Math\'ematiques, UMR 7586 du CNRS, 2, place
Jussieu, 75251 Paris Cedex 05, France}
\email{keller@math.jussieu.fr}

\date{May 10, 2007, last modified on May 14, 2008}


\begin{abstract}
  Given a triangulated $2$-Calabi-Yau category $\cc$ and a
  cluster-tilting subcategory $\ct$, the index of an object $X$ of
  $\cc$ is a certain element of the Grothendieck group of the additive
  category $\ct$. In this note, we show that a rigid object of $\cc$
  is determined by its index, that the indices of the indecomposables
  of a cluster-tilting subcategory $\ct'$ form a basis of the
  Grothendieck group of $\ct$ and that, if $\ct$ and $\ct'$ are
  related by a mutation, then the indices with respect to $\ct$ and
  $\ct'$ are related by a certain piecewise linear transformation
  introduced by Fomin and Zelevinsky in their study of cluster
  algebras with coefficients. This allows us to give a combinatorial
  construction of the indices of all rigid objects reachable from the
  given cluster-tilting subcategory $\ct$.  Conjecturally, these
  indices coincide with Fomin-Zelevinsky's $\mathbf{g}$-vectors.
\end{abstract}

\maketitle

\section{Introduction}

This note is motivated by the representation-theoretic approach to
Fomin-Zelevinsky's cluster algebras \cite{FominZelevinsky02}
\cite{FominZelevinsky03} \cite{BerensteinFominZelevinsky05}
\cite{FominZelevinsky07} developed by Marsh-Reineke-Zelevinsky
\cite{MarshReinekeZelevinsky03}, Buan-Marsh-Reineke-Reiten-Todorov
\cite{BuanMarshReinekeReitenTodorov04}, Geiss-Leclerc-Schr\"oer
\cite{GeissLeclercSchroeer05} \cite{GeissLeclercSchroeer07} and many
others, \cf \cite{BuanMarsh05} for a survey. In this approach, a
central r\^ole is played by certain triangulated $2$-Calabi-Yau
categories and by combinatorial invariants associated with their
rigid objects (we refer to \cite{IyamaReiten06}
\cite{DerksenWeymanZelevinsky07} for different approaches). Here,
our object of study is the index, which is a certain `dimension
vector' associated with each object of the given Calabi-Yau
category.

More precisely, we fix a Hom-finite $2$-Calabi-Yau triangulated category
$\cc$ with split idempotents which admits a cluster-tilting subcategory
$\ct$.
It is known from \cite{KellerReiten06} that for each object $X$ of
$\cc$, there is a triangle
\[
T_1 \to T_0 \to X \to \Sigma T_1
\]
of $\cc$, where $T_1$ and $T_0$ belong to $\ct$.
Following \cite{Palu07}, we define the index of $X$ to be the difference
$[T_0]-[T_1]$ in the split Grothendieck group $K_0(\ct)$ of the
additive category $\ct$.
We show that
\begin{itemize}
\item[-] if $X$ is rigid (\ie $\cc(X,\Sigma X)=0$), then it is
determined by its index up to isomorphism;
\item[-] the indices of the direct factors of a rigid object all lie in the
same hyperquadrant of $K_0(\ct)$ with respect to the basis given by
a system of representatives of the isomorphism classes of the indecomposables
of $\ct$;
\item[-] the indices of the direct factors of a rigid object are linearly independent;
\item[-] the indices of a system of representatives of the indecomposable
objects of any cluster-tilting subcategory $\ct'$ form a basis of
  $K_0(\ct)$. In particular, all cluster-tilting subcategories have the
  same (finite or infinite) number of pairwise non isomorphic 
  indecomposable objects.
\end{itemize}
Note that the last point was shown in Theorem I.1.8 of
\cite{BuanIyamaReitenScott07} under the additional assumption that
$\cc$ is a stable category.  We then study how the index of an object
transforms when we mutate the given cluster-tilting subcategory.  We find
that this transformation is given by the right hand side of
Conjecture~7.12 of \cite{FominZelevinsky07}, \cf section~\ref{s:h-vectors}.
This motivates the
definition of $\mathbf{g}^\dagger$-vectors as the combinatorial
counterpart to indices. If, as we expect, Conjecture~7.12 of [loc.
cit.] holds, then our $\mathbf{g}^\dagger$-vectors are identical with
the $\mathbf{g}$-vectors of [loc. cit.], whose definition we briefly
recall below.  We finally show that if
$\cc$ has a cluster-structure in the sense of
\cite{BuanIyamaReitenScott07}, then we have a bijection between
$\mathbf{g}^\dagger$-vectors and indecomposable rigid objects
reachable from $\ct$ and between $\mathbf{g}^\dagger$-clusters and
cluster-tilting subcategories reachable from $T$.

Our results are inspired by and closely related to the conjectures
of \cite{FominZelevinsky07} and the results of section~15 in
\cite{GeissLeclercSchroeer06a}. 
As a help to the reader not familiar with \cite{FominZelevinsky07},
we give a short summary of the notions introduced there which are
most relevant for us: Let $n\geq 1$ be an integer and $B$ a 
skew-symmetric integer matrix. Let $\cF$ be the field of
rational functions $\Q(x_1, \ldots, x_n, y_1, \ldots, y_n)$ in
$2n$ indeterminates. Let $\ca \subset \cF$ be the cluster algebra with
principal coefficients associated with the initial seed
$(\mathbf{x}, \mathbf{y}, B)$, where $\mathbf{x}=(x_1, \ldots, x_n)$ 
and $\mathbf{y}=(y_1, \ldots, y_n)$, \cf sections~1 and 2
of \cite{FominZelevinsky07}. As shown in Proposition~3.6 of
\cite{FominZelevinsky07}, each cluster variable of $\ca$ lies
in the ring $\Z[x_1^{\pm 1}, \ldots, x_{n}^{\pm 1}, y_1, \ldots, y_n]$.
Moreover, by Proposition~6.1 of \cite{FominZelevinsky07}, each cluster
variable of $\ca$ is homogeneous with respect to the $\Z^n$-grading
on $\Z[x_1^{\pm 1}, \ldots, x_{n}^{\pm 1}, y_1, \ldots, y_n]$ given
by 
\[
\deg(x_i)=e_i \ko \quad \deg(y_j) = -\sum_{i=1}^n b_{ij} e_i \ko
\]
where the $e_i$ form the standard basis of $\Z^n$. The $g$-vector
associated with a cluster variable $X$ is by definition the
vector $\deg(X)$ of $\Z^n$. More generally, the $g$-vector
of a cluster monomial $M$ is $\deg(M)$. Now we can state the
conjectures of \cite{FominZelevinsky07} which motivated the
above statements on the combinatorics of rigid objects:
\begin{itemize}
\item[-] different cluster monomials have different $g$-vectors
(part (1) of Conjecture~7.10 of \cite{FominZelevinsky07});
\item[-] the $g$-vectors of the variables in a fixed cluster all
lie in the same hyperquadrant of $\Z^n$ (Conjecture 6.13 of 
\cite{FominZelevinsky07});
\item[-] the $g$-vectors of the variables in a fixed cluster form a 
basis of $\Z^n$ (part (2) of Conjecture~7.10 of \cite{FominZelevinsky07});
\item[-] under a mutation of the initial cluster, the $g$-vector
of a given cluster variable transforms according to a certain
piecewise linear transformation, \cf section~\ref{s:h-vectors}
(Conjecture 7.12 of \cite{FominZelevinsky07}).
\end{itemize}
In \cite{FuKeller07}, the results of this paper have been
used to prove these conjectures for certain classes
of cluster algebras.

\section{A rigid object is determined by its index}

Let $k$ be an algebraically closed field and $\cc$ a
$\Hom$-finite $k$-linear triangulated category with
split idempotents. In particular, the decomposition
theorem holds for $\cc$: Each object decomposes into
finite sum of indecomposable objects, unique up to
isomorphism,  and indecomposable objects have local
endomorphism rings. We write $\Sigma$ for the suspension
functor of $\cc$. We suppose that $\cc$ is $2$-Calabi-Yau,
\ie that the square of the suspension functor (with its
canonical structure of triangle functor) is a Serre
functor for $\cc$. This implies that we have bifunctorial
isomorphisms
\[
D\cc(X,Y) \iso \cc\,(Y,\Sigma^2 X) \ko
\]
where $X$ and $Y$ vary in $\cc$ and $D$ denotes the
duality functor $\Hom_k(?,k)$ over the ground field.
Moreover, we suppose that $\cc$ admits a cluster-tilting
subcategory $\ct$ (called a maximal $1$-orthogonal subcategory
in \cite{Iyama05}). Recall from \cite{KellerReiten06} that this means
that $\ct$ is a full additive subcategory such that
\begin{itemize}
\item[-] $\ct$ is functorially finite in $\cc$, \ie for all
objects $X$ of $\cc$, the restrictions of the functors
$\cc(X,?)$ and $\cc(?,X)$ to $\ct$ are finitely generated, and
\item[-] an object $X$ of $\cc$ belongs to $\ct$
iff we have $\cc(T,\Sigma X)=0$ for all objects
$T$ of $\ct$.
\end{itemize}
We call an object $X$ of $\cc$ {\em rigid} if the
space $\cc(X,\Sigma X)$ vanishes.

\subsection{Rigid objects yield open orbits}
Let $X$ be a rigid object of $\cc$.
From \cite{KellerReiten07}, we know that there is a triangle
\[
\xymatrix{ T_1 \ar[r]^f & T_0 \ar[r]^{h} & X \ar[r] & \Sigma X} \ko
\]
where $T_0$ and $T_1$ belong to $\ct$.
The algebraic group $G=\Aut(T_0)\times \Aut(T_1)$ acts on
$\cc(T_1,T_0)$ via
\[
(g_0, g_1) f' = g_0 f' g_1^{-1}.
\]

\begin{lemma}
The orbit of $f$ under the action of $G$ is open in $\cc(T_1,T_0)$.
\end{lemma}

\begin{proof} It suffices to prove that the differential
of the map $g \mapsto gf$ is a surjection from $\Lie(G)$
to $\cc(T_1,T_0)$. This differential is given by
\[
(\gamma_0, \gamma_1) f = \gamma_0 f - f \gamma_1.
\]
Let $f'$ be an element of $\cc(T_1,T_0)$. Consider the following
diagram
\[
\xymatrix{
\Sigma^{-1}X \ar[r]^e &
T_1 \ar@{-->}[dl]_{\gamma_1} \ar[d]^{f'} \ar[r]^{f} &
T_0 \ar@{-->}[dl]^{\gamma_0} \ar@{-->}[d]^{\beta_0} \ar[r]^h &
X \\
T_1 \ar[r]_f & T_0 \ar[r]_h & X \ar[r] & \Sigma T_1.
}
\]
Since $X$ is rigid, the composition $hf'e$ vanishes.
So there is a $\beta_0$ such that $\beta_0 f = h f'$. Now
$h$ is a right $\ct$-approximation. So there is a
$\gamma_0$ such that $h \gamma_0 = \beta_0$. It follows
that we have
\[
h (  \gamma_0 f - f') =0.
\]
So there is a $\gamma_1$ such that
\[
\gamma_0 f  - f'= f \gamma_1.
\]
This shows that the differential of the map $g \mapsto gf$ is indeed
surjective.
\end{proof}

\subsection{Rigid objects have disjoint terms in their minimal presentations}
Let
\[
F : \cc \to \mod \ct
\]
be the functor taking an object $Y$ of $\cc$ to
the restriction of $\cc(?,Y)$ to $\ct$.
Let $X$ be a rigid object of $\cc$. Let
\[
\xymatrix{T_1 \ar[r] & T_0 \ar[r]^{h} & X \ar[r]^-{\eps} & \Sigma T_1}
\]
be a triangle such that $T_0$ and $T_1$ belong to $\ct$ and
$h$ is a minimal right $\ct$-approximation.

\begin{proposition} $T_0$ and $T_1$ do not have an indecomposable direct factor
in common.
\end{proposition}

We give two proofs of the proposition. Here is the first one:

\begin{proof} We know that
\[
FT_1 \to FT_0 \to FX \to 0
\]
is a minimal projective presentation of $FX$. Since $F$ induces
an equivalence from $\ct$ onto the category of projectives of
$\mod \ct$, it is enough to show that $FT_1$ and $FT_0$ do not
have an indecomposable factor in common. For this, it suffices
to show that no simple module $S$ occuring in the head of
$FT_0$ also occurs in the head of $FT_1$. Equivalently,
we have to show that if a simple $S$ satisfies $\Hom(FX,S)\neq 0$,
then we have $\Ext^1(FX,S)=0$. So let $S$ be a simple admitting
a surjective morphism
\[
p: FX \to S.
\]
Let $f: FT_1 \to S$ be a map representing an element in $\Ext^1(FX,S)$.
Since $FT_1$ is projective, there is a morphism $f_1 : FT_1 \to FX$
such that $p\circ f_1 = f$. Now using the fact that $F$ is essentially
surjective and full, we choose a preimage up to isomorphism
$\tilde{S}$ of $S$ and preimages
$\tilde{f}$, $\tilde{p}$ and $\tilde{f}_1$ of
$f$, $p$ and $f_1$ in $\cc$ as in the following diagram
\[
\xymatrix{
\Sigma^{-1} X \ar[r]^-{\Sigma^{-1} \eps} & T_1 \ar[dl]_{\tilde{f}_1}
\ar[d]^{\tilde{f}} \ar[r] & T_0 \ar[r] & X \\
X \ar[r]_{\tilde{p}} & \tilde{S}
}
\]
Denote by $\mod\ct$ the category of finitely presented $k$-linear functors
from $\ct^{op}$ to the category of $k$-vector spaces.
Since $F$ induces a bijection
\[
\cc(T,Y) \to (\mod \ct)(FT,FY)
\]
for all $Y$ in $\cc$, we still have $\tilde{p} \circ \tilde{f}_1 = \tilde{f}$.
The composition $\tilde{f}_1 \circ (\Sigma^{-1}\eps)$ vanishes since
we have $\cc(\Sigma^{-1} X,X)=0$. Therefore, the composition
\[
\tilde{f} \circ (\Sigma^{-1} \eps) = \tilde{p} \circ \tilde{f}_1 \circ
(\Sigma^{-1} \eps)
\]
vanishes. This implies that $\tilde{f}$ factors through the
morphism $T_1 \to T_0$. But then $f$ factors through the
morphism $FT_1 \to FT_0$ and $f$ represents $0$ in $\Ext^1(FX,S)$.
\end{proof}

Let us now give a second, more geometric, proof of the proposition:

\begin{proof}
Suppose that $T_0$ and $T_1$ have an indecomposable direct factor
$T_2$ so that we have decompositions
\[
T_0 = T_0'\oplus T_2 \mbox{ and } T_1 = T_1'\oplus T_2.
\]
For a morphism $f: T_1 \to T_0$, let
\[
\left[ \begin{array}{cc} f_{11} & f_{12} \\ f_{21} & f_{22} \end{array} \right]
\]
be the matrix corresponding to $f$ with respect to the given
decompositions. Of course, up to isomorphism, the cone on $f$ only
depends on the orbit of $f$ under the group $\Aut(T_0) \times \Aut(T_1)$.
Suppose that the cone on $f$ is isomorphic to $X$, which is rigid.
Then we know that the orbit of $f$ in $\cc(T_1,T_0)$ is open.
Hence there is some $f'$ in the orbit such that the component
$f'_{22}$ is invertible. But then, using elementary operations
on the rows and columns of the matrix of $f'$, we see that the
orbit of $f$ contains a morphism $f''$ whose matrix is diagonal
with invertible component $f''_{22}$. Clearly, the triangle
on $f''$ is not minimal. This shows that $T_1$ and $T_0$
do not have a common indecomposable factor if they are the
terms of a minimal triangle whose third term is the rigid
object $X$.
\end{proof}

\subsection{A rigid object is determined by its index}
The (split) Grothendieck group $K_0(\ct)$ of the additive category
$\ct$ is the quotient of the free group on the isomorphism
classes $[T]$ of objects $T$ of $\ct$ by the subgroup
generated by the elements of the form
\[
[T_1\oplus T_2] -[T_1] -[T_2].
\]
It is canonically isomorphic to the
free abelian group on the isomorphism classes
of the indecomposable objects of $\ct$. It contains a canonical
positive cone formed by the classes of objects of $\ct$.
Each element $c$ of $K_0(\ct)$ can be uniquely written as
\[
c=[T_0] - [T_1]
\]
where $T_0$ and $T_1$ are objects of $\ct$ without common
indecomposable factors.
Let $X$ be an object of $\cc$. Recall that its index \cite{Palu07} is the
element
\[
\ind(X)=[T_0] - [T_1]
\]
of $K_0(\ct)$ where $T_0$ and $T_1$ are objects of $\ct$ which
occur in an arbitrary triangle
\[
T_1 \to T_0 \to X \to \Sigma T_1.
\]
Now suppose that $X$ is rigid. We know that if we choose the
above triangle minimal, then $T_0$ and $T_1$ do not have common indecomposable
factors. Thus they are determined by $\ind(X)$. Moreover, since the
$\cc(T_1,T_0)$ is an irreducible variety (like any finite-dimensional
vector space),
each morphism $f: T_1 \to T_0$ whose orbit under the
group $\Aut(T_0)\times \Aut(T_1)$ is open yields a cone isomorphic
to $X$. Thus up to isomorphism, $X$ is determined by $\ind(X)$.
In fact, $X$ is
isomorphic to the cone on a general morphism $f: T_1 \to T_0$ between
the objects $T_0$ and $T_1$ without a common indecomposable
factor such that $\ind(X)=[T_0]-[T_1]$. We have proved the

\begin{theorem} \label{thm:injection} The map $X \mapsto \ind(X)$ induces an injection
from the set of isomorphism classes of rigid objects of $\cc$ into
the set $K_0(\ct)$.
\end{theorem}

\noindent
This theorem was inspired by part~(1) of conjecture~7.10 in
\cite{FominZelevinsky07}.

\subsection{Direct factors of rigid objects have sign-coherent indices}
Let $A$ be a free abelian group endowed with a basis $e_i$, $i\in I$.
A subset $X \subset A$ is {\em sign-coherent} if, for all elements
$x,y\in X$ and for all $i\in I$, the sign of the component $x_i$ in the decomposition
\[
x=\sum x_i e_i
\]
agrees with the sign of $y_i$, \cf Definition~6.12 of \cite{FominZelevinsky07}.
This means that the set $X$ is entirely contained in a hyperquadrant
of $A$ with respect to the given basis $e_i$, $i\in I$.
Now consider the free abelian
group $K_0(\ct)$ endowed with the basis formed by the classes of
indecomposable objects of $\ct$. Suppose that $X$ is a rigid object
of $\cc$. We claim that the set of indices of the direct factors of $X$ is
sign-coherent. Indeed, let $U$ and $V$ be direct factors of $X$. Choose
minimal triangles
\[
T_1^U \to T_0^U \to U \to \Sigma T_1^U \mbox{ and }
T_1^V \to T_0^V \to V \to \Sigma T_1^V \ko
\]
where the $T_i^U$ and $T_i^V$ belong to $\ct$. Then the triangle
\[
T_1^U \oplus T_1^V \to T_0^U\oplus T_0^V \to U\oplus V
\to \Sigma(T_1^U \oplus T_1^V)
\]
is minimal. Since $U\oplus V$ is rigid, the two terms
$T_1^U \oplus T_1^V$  and  $T_0^U\oplus T_0^V$ do not have indecomposable
direct factors in common. In particular, whenever an indecomposable
object occurs in $T_0^U$ (resp. $T_1^U$), it does not occur in
$T_1^V$ (resp. $T_0^V$). This shows that $\ind(U)$ and $\ind(V)$
are sign-coherent. This property is to be compared with
conjecture~6.13 of \cite{FominZelevinsky07}.

\subsection{Indices of factors of rigid objects are linearly independent}
Let $X$ be a rigid object of $\cc$ and let $X_i$, $i\in I$, be a
finite family of indecomposable direct factors of $X$ which
are pairwise non isomorphic. We claim that the
elements $\ind(X_i)$, $i\in I$, are linearly independent in $K_0(\ct)$.
Indeed, suppose that we have a relation
\[
\sum_{i\in I_1} c_i \ind(X_i) = \sum_{j\in I_2} c_j \ind(X_j)
\]
for two disjoint subsets $I_1$ and $I_2$ of $I$ and positive
integers $c_i$ and $c_j$. Then the rigid
objects
\[
\bigoplus_{i\in I_1} X_i^{c_i} \mbox{ and } \bigoplus_{j\in I_2} X_j^{c_j}
\]
have equal indices. So they are isomorphic. Since $I_1$ and $I_2$
are disjoint, all the $c_i$ and $c_j$ have to vanish.

\subsection{The indices of the indecomposables of a cluster tilting subcategory
form a basis} The following theorem was inspired by
part~(2) of conjecture~7.10 of \cite{FominZelevinsky07}.

\begin{theorem} \label{thm:basis}
Let $\ct'$ be another tilting subcategory of $\cc$. Then the
elements $\ind(T')$, where $T'$ runs through a system of representatives
of the isomorphism classes of
indecomposables of $\ct'$, form a basis of the free abelian
group $K_0(\ct)$.
\end{theorem}

\begin{proof}
Indeed, we already know that the $\ind(T')$ are
linearly independent. So it is enough to show that the subgroup they
generate contains $\ind(T)$ for each indecomposable $T$ of $\ct$.
Indeed, let $T$ be an indecomposable of $\ct$ and let
\[
T \to T'_1 \to T'_0 \to \Sigma T
\]
be a triangle with $T'_i$ in $\ct'$ (this triangle allows to
compute the index of $\Sigma T$ with respect to $\ct'$). Then the map
$FT'_1 \to FT'_0$ is surjective and therefore, we have
\[
\ind(T) - \ind(T'_1) + \ind(T'_0) = 0
\]
by Proposition~6 of \cite{Palu07}. Thus, $\ind(T)$ is in the subgroup
of $K_0(\ct)$ generated by the $\ind(T')$, where $T'$ runs through
the indecomposables of $\ct'$.
\end{proof}

\section{How the index transforms under change of cluster-tilting subcategory}

Let $\ct'$ be another cluster-tilting subcategory. Suppose
that $\ct$ and $\ct'$ are {\em related by a mutation}, \ie
there is an indecomposable $S$ of $\ct$ and an indecomposable
$S^*$ of $\ct'$ such that, if $\indec$ denotes the set of
isomorphism classes of indecomposables, we have
\[
\indec(\ct')=\indec(\ct) \setminus \{S\} \cup \{S^*\},
\]
and that there exist triangles
\[
S^* \to B \to S \to \Sigma S^* \mbox{ and }
S \to B' \to S^* \to \Sigma S
\]
with $B$ and $B'$ belonging to $\ct\cap\ct'$,
\cf \eg \cite{BuanMarshReinekeReitenTodorov04} \cite{GeissLeclercSchroeer07}
\cite{IyamaYoshino06}.
We define two linear maps
\[
\phi_+ : K_0(\ct) \to K_0(\ct') \mbox{ and }
\phi_- : K_0(\ct) \to K_0(\ct')
\]
which both send each indecomposable $T''$ belonging to
both $\ct$ and $\ct'$ to itself and such that
\[
\phi_+(S) = [B] - [S^*] \mbox{ and } \phi_-(S) = [B'] - [S^*].
\]
For an object $X$ of $\cc$, we denote by $\ind_\ct(X)$ the
index of $X$ with respect to $\ct$ and by $[\ind_\ct(X):S]$
the coefficient of $S$ in the decomposition of $\ind_\ct(X)$
with respect to the basis given by the indecomposables of $\ct$.
The following theorem is inspired by Conjecture~7.12 of
\cite{FominZelevinsky07}.

\begin{theorem} \label{thm:index-transform}
Let $X$ be a rigid object of $\cc$. We have
\[
\ind_{\ct'}(X) = \left\{ \begin{array}{ll}
\phi_+(\ind_\ct(X)) &  \mbox{if}\quad [\ind_\ct(X):S] \geq 0 \;\; ;\\
\phi_-(\ind_\ct(X)) &  \mbox{if}\quad [\ind_\ct(X):S] \leq 0 .
\end{array} \right.
\]
\end{theorem}

\begin{proof}
Let
\[
T_1 \to T_0 \to X \to \Sigma T_1
\]
be a triangle with $T_0$ and $T_1$ in $\ct$.
Suppose first that $S$ occurs neither as a
direct factor of $T_1$ nor of $T_0$. Then clearly the
triangle yields both the index of $X$ with respect
to $\ct$ and with respect to $\ct'$ and we have
\[
\phi_+(\ind_\ct(X)) = \phi_-(\ind_{\ct}(X)) = \ind_{\ct'}(X).
\]
Now suppose that the multiplicity $[\ind_\ct(X):S]$
equals a positive integer $i\geq 1$. This means that $S$
occurs with multiplicity $i$ in $T_0$ but does not occur
as a direct factor of $T_1$. Choose a decomposition
$T_0 = T_0'' \oplus S^i$. From the octahedron constructed over
the composition
\[
T_0''\oplus B^i \to T_0''\oplus S^i \to X \ko
\]
we extract the following commutative diagram, whose rows
and columns are triangles
\[
\xymatrix{
\Sigma S^{*i} \ar[r]^{\id} & \Sigma S^{*i} & & \\
T_1 \ar[u] \ar[r] & T_0''\oplus S^i \ar[u] \ar[r] & X \ar[r] & \Sigma T_1 \\
T_1' \ar[u] \ar[r] & T_0''\oplus B^i \ar[r] \ar[u] & X \ar[u]_{\id} \ar[r] &
                                                       \Sigma T_1' \ar[u] \\
S^{*i} \ar[u] \ar[r]^{\id} & S^{*i}. \ar[u] & &
}
\]
Since there are no non zero morphisms from $T_1$ to $\Sigma S^{*i}$
($T_1$ and $S^*$ belong to $\ct'$), the leftmost column is a split triangle
and $T_1'$ is isomorphic to $S^{*i}\oplus T_1$. Thus, the third
line yields the index of $X$ with respect to $\ct'$, which
equals
\[
\ind_{\ct'}(X) = [T_0''\oplus B^i] - [T_1'] = [T_0''] - [T_1] + i([B]-[S^*])
= \phi_+(\ind_\ct(X)).
\]
Finally, suppose that the multiplicity $[\ind_\ct(X):S]$ is equals
a negative integer $-i\leq -1$. This means that $S$ occurs with
multiplicity $i$ in $T_1$ but does not occur in $T_0$.
Choose a decomposition $T_1=T_1'' \oplus S^i$. From
the octahedron over the composition
\[
\Sigma^{-1} X \to T_1'' \oplus S^i \to T_1''\oplus B'^i \ko
\]
we extract the following diagram, whose rows and columns are
triangles
\[
\xymatrix{
                                & \Sigma^{-1} S^{*i} \ar[r]^{\id} \ar[d]  & \Sigma^{-1} S^{*i}  \ar[d] &                \\
\Sigma^{-1} X \ar[d]_\id \ar[r] & T_1''\oplus S^i \ar[r] \ar[d]           & T_0 \ar[r] \ar[d]          & X \ar[d]^\id   \\
\Sigma^{-1} X \ar[r]            & T_1''\oplus B'^i \ar[r] \ar[d]          & T_0' \ar[d] \ar[r]         & X              \\
                                & S^{*i} \ar[r]_\id                       & S^{*i} &
}
\]
Since there are no non zero morphisms from $\Sigma^{-1} S^{*i}$
to $T_0$ ($S^*$ and $T_0$ belong to $\ct'$), the object $T_0'$ is
isomorphic to $T_0\oplus S^i$ and we can read $\ind_{\ct'}(X)$
off the third line of the diagram:
\[
\ind_{\ct'}(X) = [T_0'] - [T_1''\oplus B'^i] = [T_0\oplus S^{*i}] - [T_1''] - i [B']
= [T_0] - [T_1''] -i ( [B'] - [S^*]) = \phi_-(\ind_\ct(X)).
\]

\end{proof}

\section{$\mathbf{g}^\dagger$-vectors and $\mathbf{g}^\dagger$-clusters}
\label{s:h-vectors}

In this section, we recall fundamental constructions from
\cite{FominZelevinsky07} in a language adapted to our
applications. We will define $\mathbf{g}^\dagger$-vectors using
the right hand side of Conjecture~7.12 of [loc. cit.].
If, as we expect, this conjecture holds, then our $\mathbf{g}^\dagger$-vectors
are identical with the $\mathbf{g}$-vectors of [loc. cit.].

Let $Q$ be a quiver. Thus $Q$ is given by a set of vertices $I=Q_0$,
a set of arrows $Q_1$ and two maps $s$ and $t$ from $Q_1$ to $I=Q_0$
taking an arrow to its source, respectively its target. We assume that
$Q$ is {\em locally finite}, \ie for each given vertex $i$ of $Q$ there
are only finitely many arrows $\alpha$ such that $s(\alpha)=i$ or $t(\alpha)=i$.
Moreover, we assume that
$Q$ has no loops (\ie arrows $\alpha$ such that $s(\alpha)=t(\alpha)$)
and no $2$-cycles (\ie pairs of distinct arrows $\alpha\neq \beta$
such that $s(\alpha)=t(\beta)$ and $t(\beta)=s(\alpha)$). The
quiver $Q$ is thus determined by the set $I$ and the skew-symmetric
integer matrix $B=(b_{ij})_{I\times I}$ such that, whenever the coefficient $b_{ij}$
is positive, it equals the number of arrows from $i$ to $j$ in $Q$.
Notice that if, for an integer $x$, we write $[x]_+=\max(x,0)$,
then the number of arrows from $i$ to $j$  in $Q$ is $[b_{ij}]_+$.
The {\em mutation} $\mu_k(Q)$ of $Q$ at a vertex $k$ is
by definition the quiver with vertex set $I$ whose numbers
of arrows are given by the mutated matrix $B'=\mu_k(B)$ as defined,
for example, in definition~2.4 of \cite{FominZelevinsky07}:
\[
b'_{ij} = \left\{ \begin{array}{ll}
-b_{ij} & \mbox{if $i=k$ or $j=k$;} \\
b_{ij}+\sgn(b_{ik})[b_{ij} b_{kj}]_+ & \mbox{otherwise.}
\end{array} \right.
\]
As in definition~2.8 of \cite{FominZelevinsky07}, we let
$\mathbb{T}=\mathbb{T}_I$ be the regular tree whose edges are labeled
by the elements of $I$ such that for each vertex $t$ and each
element $k$ of $I$, there is precisely one edge incident with
$t$ and labeled by $k$. We fix a vertex $t_0$ of $\mathbb{T}$
and define $Q_{t_0}=Q$. Clearly, there is a unique map assigning
a quiver $Q_t$ to each vertex $t$ such that if $t$ and $t'$
are linked by an edge labeled by $k$, we have $Q_{t'}=\mu_k(Q_t)$.
In analogy with the terminology of \cite{FominZelevinsky07},
we call the map $t \mapsto Q_t$ the {\em quiver pattern} associated
with $t_0$ and $Q$.

Now for each vertex $t$ of $\mathbb{T}$, we define $K_t$ to
be the free abelian group on the symbols $e_i^t$, $i\in I$.
For two vertices $t$ and $t'$ linked by an edge labeled $k$,
we let
\[
\phi^+_{t',t} : K_t \to K_{t'} \mbox{ respectively }
\phi^-_{t',t} : K_t \to K_{t'}
\]
be the linear map sending $e_j^t$ to $e_j^{t'}$ for
each $j\neq k$ and sending $e_k$ to
\[
-e_k^{t'}+ \Sigma_j [b^t_{jk}]_+ \, e_j^{t'} \mbox{ respectively }
-e_k^{t'} + \Sigma_j [b^t_{kj}]_+ \, e_j^{t'} \ko
\]
where $(b^t_{ij})$ is the skew-symmetric matrix associated
with the quiver $Q_t$. We define the piecewise linear transformation
\[
\phi_{t',t} : K_t \to K_{t'}
\]
to be the map whose restriction to the halfspace of
elements with positive $e_k^t$-coordinate is $\phi^+_{t',t}$
and whose restriction to the opposite halfspace
is $\phi^-_{t',t}$.
Thus, the image of an element $g$ with coordinates
$g_j$, $j\in I$, is the element $g'$ with coordinates
\[
g'_j = \left\{ \begin{array}{ll} -g_j & \mbox{if $j= k$;} \\
g_j+[b^t_{jk}]_+\, g_k & \mbox{if $j\neq k$ and $g_k\geq 0$;} \\
g_j+[b^t_{kj}]_+\, g_k & \mbox{if $j\neq k$ and $g_k\leq 0$.}
\end{array} \right.
\]
It is easy to check that this rule agrees with formula~(7.18)
in Conjecture~7.12 of \cite{FominZelevinsky07}.

If $t$ and $t'$ are two arbitrary vertices of $\mathbb{T}$,
there is a unique path
\[
\xymatrix{t=t_1 \ar@{-}[r] & t_2 \ar@{-}[r] & \ldots & t_N=t' \ar@{-}[l] }
\]
of edges leading from $t$ to $t'$ and we define $\phi_{t',t}$ to
be the composition
\[
\phi_{t_N,t_{N-1}} \circ \cdots  \circ \phi_{t_2,t_1}.
\]
For a vertex $t$ of $\mathbb{T}$ and a vertex $l$ of $Q$, the
{\em $\mathbf{g}^\dagger$-vector $\mathbf{g}^\dagger_{l,t}$} is the element of the abelian
group $K_{t_0}$ defined by
\[
\mathbf{g}^\dagger_{l,t} = \phi_{t_0, t}(e^t_l).
\]
The {\em $\mathbf{g}^\dagger$-cluster} associated with a vertex $t$ of $\mathbb{T}$
is the set of $\mathbf{g}^\dagger$-vectors $\mathbf{g}^\dagger_{l,t}$, $l\in I$.
If Conjecture~7.12 of \cite{FominZelevinsky07} holds for the cluster
algebra with principal coefficients associated with the matrix $B$,
then it is clear that in the notations of formula~(6.4) of \cite{FominZelevinsky07},
we have
\[
\mathbf{g}^\dagger_{l,t} = \mathbf{g}_{l,t}
\]
for all vertices $t$ of $\mathbb{T}$ and all $l\in I$, , \ie the $\mathbf{g}^\dagger$-vectors
equal the $\mathbf{g}$-vectors for the cluster algebra with principal
coefficients associated with the skew-symmetric matrix $B$.

\section{Rigid objects in $2$-Calabi-Yau categories with cluster structure}
\label{s:rigid-objects}

Let $\cc$ be a $\Hom$-finite $2$-Calabi-Yau category with a cluster-tilting
subcategory $\ct$. Let $Q=Q(\ct)$ be the {\em quiver of $\ct$}. Recall that this
means that the vertices of $Q$ are the isomorphism classes of indecomposable
objects of $\ct$ and that the number of arrows from the isoclass of
$T_1$ to that of $T_2$ equals the dimension of the space of irreducible
morphisms
\[
\mathsf{irr}(T_1, T_2) = \mathsf{rad}(T_1, T_2)/\mathsf{rad}^2(T_1,T_2) \ko
\]
where $\mathsf{rad}$ denotes the radical of $\ct$, \ie the ideal such
that $\mathsf{rad}(T_1,T_2)$ is formed by all non isomorphisms from $T_1$ to $T_2$.

We make the following {\em assumption on $\cc$}: For each cluster-tilting
subcategory $\ct'$ of $\cc$, the quiver $Q(\ct')$ does not have loops
or $2$-cycles. We refer to section~1, page~11 of \cite{BuanIyamaReitenScott07} for
a list of classes of examples where this assumption holds.
By theorem~1.6 of \cite{BuanIyamaReitenScott07}, the assumption
implies that the cluster-tilting subcategories of $\cc$
determine a cluster structure for $\cc$. Let us
recall what this means:
\begin{itemize}
\item[1)] For each cluster-tilting subcategory $\ct'$ of
$\cc$ and each indecomposable $S$ of $\ct'$, there is a unique (up
to isomorphism) indecomposable $S^*$ not isomorphic to $M$ and such
that the additive subcategory $\ct''=\mu_S(\ct')$ of $\cc$ with
\[
\mathsf{indec}(\ct'')=\mathsf{indec}(\ct')\setminus\{S\} \cup \{S^*\}
\]
is a cluster-tilting subcategory;
\item[2)] the space of morphisms
from $S$ to $\Sigma S^*$ is one-dimensional and in the non-split
triangles
\[
S^* \to B \to S \to \Sigma S^* \mbox{ and } S \to B' \to S^* \to \Sigma S
\]
the objects $B$ and $B'$ belong to $\ct'\cap\ct''$;
\item[3)] the multiplicity
of an indecomposable $L$ of $\ct'\cap\ct''$ in $B$ equals the number of arrows from
$L$ to $S$ in $Q(\ct')$  and that from $S^*$ to $L$ in $Q(\ct'')$;
the multiplicity of $L$ in $B'$ equals the
number of arrows from $S$ to $L$ in $Q(\ct')$ and that from
$L$ to $S^*$ in $Q(\ct'')$;
\item[4)] finally, we have $Q(\ct'')=\mu_S(Q(\ct'))$.
\end{itemize}

Let $Q=Q(\ct)$ be the quiver of $\ct$. Notice that its set of
vertices is the set $Q_0=I$ of isomorphism classes of indecomposables
of $\ct$. Let $\mathbb{T}$ be the regular tree associated with $Q$
as in section~\ref{s:h-vectors}. We fix a vertex $t_0$ of $\mathbb{T}$
and put $\ct_{t_0}=\ct$. For two cluster tilting subcategories
$\ct'$ and $\ct''$ as above, let $\psi_{\ct'',\ct'}: \indec(\ct') \to \indec(\ct'')$
be the bijection taking $S$ to $S^*$ and fixing all other indecomposables.

Thanks to point 1), with each
vertex $t$ of $\mathbb{T}$, we can associate
\begin{itemize}
\item[a)] a unique cluster-tilting
subcatgory $\ct_t$ and
\item[b)] a unique bijection
\[
\psi_{t,t_0}: \indec(\ct_{t_0}) \to \indec(\ct_t)
\]
\end{itemize}
such that $\ct_{t_0}=\ct$ and that, whenever
two vertices $t$ and $t'$ are linked by an edge labeled by
an indecomposable $S$ of $\ct=\ct_{t_0}$, we have
\begin{itemize}
\item[a)] $\ct_{t'}=\mu_{S'}(\ct_{t})$, where $S'=\psi_{t,t_0}(S)$, and
\item[b)] $\psi_{t',t_0}=\psi_{t',t}\circ \psi_{t,t_0}$.
\end{itemize}
Moreover, thanks to point 4), the map $t\mapsto Q(\ct_t)$
is the quiver-pattern associated with $Q$ 
and $t_0$ in section~\ref{s:h-vectors}.
Notice that the group $K_0(\ct)$ with the basis formed by
the isomorphism classes of indecomposables canonically
identifies with the free abelian group $K_{t_0}$
of section~\ref{s:h-vectors}. We define a cluster-tilting
subcategory $\ct'$ to be {\em reachable from $\ct$} if
we have $\ct'=\ct_t$ for some vertex $t$ of the tree $\mathbb{T}$.
We define a rigid indecomposable $M$ to be {\em reachable from
$\ct$} if it belongs to a cluster-tilting subcategory which
is reachable from $\ct$.

\begin{theorem} \begin{itemize}
\item[a)] The index $\ind(M)$ of a rigid indecomposable reachable
from $\ct$ is a $\mathbf{g}^\dagger$-vector 
and the map $M\mapsto \ind(M)$ induces a
bijection from the set of isomorphism classes of rigid indecomposables
reachable from $\ct$ onto the set of $\mathbf{g}^\dagger$-vectors.
\item[b)] Under the bijection $M\mapsto \ind(M)$
of a), the cluster-tilting subcategories reachable from $\ct$
are mapped bijectively to the $\mathbf{g}^\dagger$-clusters.
\end{itemize}
\end{theorem}

\begin{proof} a) By assumption, there is a vertex $t$ of $\mathbb{T}$
such that $M$ belongs to $\ct_t$. Now we use theorem~\ref{thm:index-transform}
and induction on the length of the path joining $t_0$ to $t$ in the
tree $\mathbb{T}$ to conclude that
\[
\ind(M)=\mathbf{g}^\dagger_{M',t}\ko \mbox{ where } M=\psi_{t,t_0}(M').
\]
This formula  shows that the map $M \mapsto \ind(M)$ is
a well-defined surjection onto the set of $\mathbf{g}^\dagger$-vectors.
By theorem~\ref{thm:injection}, the
map $M\mapsto \ind(M)$ is also injective.
b) By assumption, a reachable cluster-tilting subcategory
$\ct'$ is of the form $\ct_t$ for some vertex $t$ of
the tree $\mathbb{T}$. Thus its image is the $\mathbf{g}^\dagger$-cluster
associated with $t$. This shows that the map is well-defined
and surjective. It follows from a) that it is also injective.
\end{proof}

\section*{Acknowledgments}

The second-named author would like to thank Jan Schr\"oer for
stimulating discussions. Both authors are grateful to Andrei Zelevinsky
for helpful comments on a previous version of this article.

\def\cprime{$'$}
\providecommand{\bysame}{\leavevmode\hbox
to3em{\hrulefill}\thinspace}
\providecommand{\MR}{\relax\ifhmode\unskip\space\fi MR }
\providecommand{\MRhref}[2]{%
  \href{http://www.ams.org/mathscinet-getitem?mr=#1}{#2}
} \providecommand{\href}[2]{#2}

\end{document}